\documentclass[11pt]{article}
\usepackage{amssymb,amsfonts,latexsym,amsthm,graphics}

\mathsurround.8pt\newcommand{\n}{\noindent}
\newcommand\pr{\sb\n\textit{Proof. }}  
\newcommand\ex{\mb\n\textit{Example. }}
\newcommand\re{\mb\n\textit{Remark. }}
\newcommand{\tr}{\hbox{tr}}\renewcommand{\Re}{\hbox{Re}}
\renewcommand{\a}{\alpha}\newcommand{\oa}{\ol\a}
\renewcommand{\b}{\beta}
\newcommand{\la}{\lambda}\newcommand{\W}{\Omega}
\newcommand{\C}{\mathbb{C}}\newcommand{\R}{\mathbb{R}}
\newcommand{\bD}{\mathbb{D}}\newcommand{\bM}{\mathbb{M}}
\newcommand{\ga}{\gamma}\newcommand{\w}{\omega} 
\newcommand{\si}{\sigma}\newcommand{\ow}{\ol\omega}
\newcommand\hg{\widehat{\gamma}}
\newcommand{\g}{\mathfrak{g}}\newcommand{\G}{G}
\newcommand{\X}{X}\newcommand{\XoX}{X\ol X}
\newcommand{\La}{\Lambda}\newcommand{\AC}{\mathcal{A}}
\newcommand{\cC}{\mathcal{C}}\newcommand{\cO}{\mathcal{O}}
\newcommand{\cU}{\mathcal{U}}\newcommand{\cV}{\mathcal{V}}
\newcommand{\cS}{\mathcal{S}}\newcommand{\pZ}{\mathcal{Z}'_-}
\newcommand{\es}{\emptyset}\newcommand{\sm}{\setminus}
\newcommand{\cCg}{\cC(\g)}\newcommand{\cCpg}{\cC^+(\g)}
\newcommand{\cCgg}{\mathcal{C}^+(\g,g)}\newcommand{\ret}{r}
\newcommand\bu{{\hbox{\footnotesize$\bullet$}}}
\newcommand{\cCf}{\cC^\bu}\newcommand{\cCi}{\cC^\infty}
\newcommand{\ACf}{\AC^\bu}\newcommand{\p}{p}
\newcommand\hJ{\widehat{J}}
\newtheorem{theorem}{Theorem}[section] 
\newtheorem{lemma}[theorem]{Lemma}
\newtheorem{definition}[theorem]{Definition}
\newtheorem{corollary}[theorem]{Corollary}
\newtheorem{proposition}[theorem]{Proposition}
\def\E{\raise1pt\hbox{$\bigwedge$}}\newcommand{\we}{\wedge}
\renewcommand\qed{\hfill$\Box$\medbreak} \def\qq{\qquad}\def\q{\quad}
\def\mb{\medbreak}\def\sb{\smallbreak}
\newcommand{\al}{\langle}\newcommand{\ar}{\rangle}
\newcommand{\ol}{\overline}\newcommand{\op}{\oplus}
\newcommand{\be}{\begin{equation}}\newcommand{\ee}{\end{equation}}
\newcommand{\ba}{\begin{array}}\newcommand{\ea}{\end{array}}

\renewcommand{\ge}{\geqslant}\renewcommand{\le}{\leqslant}

\begin{document}\parskip1pt

\title{\Large\bf Complex structures on the Iwasawa manifold}

\author{Georgios Ketsetzis$^1$ \and Simon Salamon}\date{}\maketitle

\begin{abstract}\n We identify the space of left-invariant oriented 
complex structures on the complex Heisenberg group, and prove that it has the
homotopy type of the disjoint union of a point and a 2-sphere.
\end{abstract}

\footnotetext[1]{\scriptsize supported by the Alexander S.~Onassis
Public Benefit Foundation and State Scholarships Foundation, Greece}

\subsection*{Introduction}

It is well known that every even-dimensional compact Lie group has a
left-invariant complex structure \cite{Sam,Wang}. By contrast, not all
nilpotent groups admit left-invariant complex structures. In 6 real dimensions
there are 34 isomorphism classes of simply-connected nilpotent Lie groups, and
the study~\cite{S6} reveals that 18 of these admit invariant complex
structures. The complex Heisenberg group $\G$ possesses a particularly rich
structure in this regard, since it has a 2-sphere of abelian complex structures
in addition to its standard bi-invariant complex structure $J_0$.

The Iwasawa manifold $\bM=\Gamma\backslash\G$ is a compact quotient of $\G$,
and any left-invariant tensor on $\G$ induces a tensor on $\bM$. As explained
in \S2, studies of Dolbeault cohomology suggest that the moduli space of
complex structures on $\bM$ is determined by the space of left-invariant
complex structures on $\G$. The set of such structures compatible with a
standard metric $g$ and orientation is the union of $\{J_0\}$ and the 2-sphere
already mentioned \cite{AGS}. The present paper shows that this description
remains valid at the level of homotopy when one no longer insists on
compatibility with $g$. This requires a new approach, in which complex
structures are described by a basis of $(1,0)$-forms in echelon form (see
Proposition 2.3). Similar techniques can be applied to other Lie groups and
nilmanifolds, though we refer the reader to \cite{PP} for related studies.

We work mainly with the Lie algebra $\g$ of $\G$, and regard left-invariant
differential forms on $\G$ as elements of $\E^k\g^*$. A special feature of the
space $\cCg$ of all invariant complex structures on $\bM$ is that any $J$ in
$\cCg$ is compatible with the fibration of $\bM$ as a $T^2$ bundle over
$T^4$. Algebraically, this amounts to asserting that the 4-dimensional kernel
$\bD$ of $d:\g^*\to\E^2\g^*$ is necessarily $J$-invariant. As we show in \S2,
the essential features of an invariant complex structure $J$ are captured by
its restriction to $\bD$, and are encoded in a complex $2\times2$ matrix
$X$. In this way, topological questions are related to properties of the
eigenvalues of $\XoX$ and some matrix analysis described in \cite{HH}.

The orientation of the restriction of an almost complex structure $J$ to $\bD$
determines two connected components of $\cCg$ that we study separately. We
establish global complex coordinates on the component $\cC_+$ containing the
complex structure induced by $J_0$, and show that it has the structure of a
contractible complex 6-dimensional manifold. By exploiting an $SU(2)$ action on
the second component $\cC_-$, we prove that this retracts onto the 2-sphere of
negatively-oriented orthogonal almost complex structures on $\bD$.

\setcounter{theorem}0
\subsection{Preliminaries}

The Iwasawa manifold $\bM$ is defined as the quotient $\Gamma\backslash\G$,
where \[\G=\left\{\left(\ba{ccc}1&z_1&z_3\\0&1& z_2\\0 & 0 &1\ea\right):
z_i\in\C\right\}\] is the complex Heisenberg group and $\Gamma$ is the lattice
defined by taking $z_1,z_2,z_3$ to be Gaussian integers, acting by left
multiplication. We shall regard $\bM$ as a real manifold of dimension $6$,
and we let $\g$ denote the real 6-dimensional Lie algebra associated to $\G$.

An \textit{invariant} complex structure on $\bM$ is by definition one induced
from a left-invariant complex structure on the real Lie group underlying
$\G$. Such a structure is invariant by the action of the centre $Z$ of $\G$,
that persists on $\bM$ ($Z$ consists of matrices for which $z_1=0=z_2$). The
set of such structures can be identified with the set $\cCg$ of almost complex
structures on the real Lie algebra $\g$ that satisfy the Lie algebraic
counterpart \[[JX,JY]=[X,Y]+ J[JX,Y]+J[X,JY]\] of the Newlander-Nirenberg
integrability condition.

The natural complex structure $J_0$ of $\G$, for which $z_1,z_2,z_3$ are
holomorphic, is a point of $\cCg$ that satisfies the stronger condition
$[JX,Y]=J[X,Y]$. It induces a bi-invariant complex structure of $\G$ that
therefore passes to a $G$-invariant complex structure on $\bM$. We shall denote
by $\cCpg$ the subset consisting of complex structures inducing the same
orientation as $J_0$.

The 1-forms \be\label{zz}\w^1=dz_1,\q\w^2=dz_2,\q\w^3=-dz_3+z_1dz_2,\ee are
left-invariant on $\G$. Define a basis $\{e^1,\ldots,e^6\}$ of \textit{real}
1-forms by setting \be\label{ch4-2}\w^1=e^1 + ie^2, \q\w^2=e^3+
ie^4,\q\w^3=e^5+ ie^6.\ee These 1-forms are pullbacks of corresponding 1-forms
on the quotient $\bM$, which we denote by the same symbols. They satisfy
\be\label{gs}\left\{\ba{ll} de^i=0,\qq 1\le i\le 4,\\[1mm]
de^{5}=e^{13}+e^{42},\\[1mm] de^{6}=e^{14}+e^{23}.\ea\right.\ee Here, we make
use of the notation $e^{ij}=e^i\we e^j$.

Let $T^k\cong\R^k/\mathbb{Z}^k$ denote a real $k$-dimensional torus.  Then
$\bM$ is the total space of a principal $T^2$-bundle over $T^4$. The mapping
$\p\colon\bM\to T^4$ is induced from $(z_1,z_2,z_3)\mapsto(z_1,z_2)$. The
space of invariant 1-forms annihilating the fibres of $p$ is \[\bD=\al
e^1,e^2,e^3,e^4\ar=\ker(d:\g^*\to\E^2\g^*),\] and this 4-dimensional subspace 
of $\g^*$ will play a crucial role in the theory.

\begin{theorem}\label{JD} Let $J$ be any invariant complex structure 
$J$ on $\bM$. Then $\p$ induces a complex structure $\hJ$ on $T^4$ such that
$\p\colon(\bM,J)\to(T^4,\hJ)$ is holomorphic.\end{theorem}

\pr Let $J$ be an element of $\cC(\g)$. The essential point is that $\bD$ is
$J$-invariant. Once this is established, it suffices to define $\hJ$ to be
the $T^4$-invariant complex structure determined on cotangent vectors by
$J|_\bD$. The pullback of a $(1,0)$-form on $T^4$ is then an invariant
$(1,0)$-form on $\bM$.

Let $\La$ denote the space of $(1,0)$-forms relative to $J$.  Then \[\dim(\al
e^1,e^2,e^3,e^4,e^5\ar_c\cap\La)=2.\] If $\dim(\bD_c\cap\La)=2$ then
$J\bD=\bD$, as required. If not, there exists a $(1,0)$-form $\delta+e^5$ with
$\delta\in\bD$. This implies that \[ de^5\in\La^{2,0}\op\La^{1,1},\] and
consequently that $de^5\in\La^{1,1}$.  Similarly for $e^6$, and thus
\[\w^1\we\w^2=d\w^3=de^5+ide^6\in\La^{1,1},\] implying that $J\w^1\we
J\w^2=\w^1\we\w^2$ and hence \[\al J\w^1,J\w^2\ar=\al\w^1,\w^2\ar.\] Thus, the
subspace $\al\w^1,\w^2\ar$ is $J$-invariant, and $J\bD=\bD$.\qed

Decreeing the 1-forms $e^i$ to be orthonormal determines a left-invariant
metric \be\label{g} g=\sum_{i=1}^6 e^i\otimes e^i\ee on $\G$. This induces
metrics on $T^4$ and $\bM$ (that we also denote by $g$) for which $\p$ is a
Riemannian submersion. The subset $\cCgg$ of $\cCpg$ corresponding to
$g$-orthogonal oriented complex structures is now easy to describe in terms of
Theorem~\ref{JD}.

\begin{lemma} The restriction of the mapping $J\mapsto\hJ$ to $\cCgg$ is 
injective.\end{lemma}

\pr We need to describe the set of invariant orthogonal complex structures on
$T^4$ in terms of 2-forms on $\bD$. First recall that an element of $\cCgg$ is
determined by the corresponding \textit{fundamental $2$-form} $\ga$ satisfying
$\ga(X,Y)=g(JX,Y)$. Given $J$ in $\cCgg$, both $\bD$ and $\bD^\perp=\al
e^5,e^6\ar$ are $J$-invariant and there exists an orthonormal basis
$\{f^1,f^2,Jf^1,Jf^2\}$ of $\bD$ for which \be\label{56}\ga=f^1\we Jf^1+f^2\we
Jf^2\pm e^5\we e^6.\ee Then the fundamental 2-form of $\hJ$ is
\be\label{hw}\hg=f^1\we Jf^1+f^2\we Jf^2.\ee The fact that the overall
orientation of $J$ on $\g$ is positive then determines uniquely the sign in
(\ref{56}).\qed

To continue the discussion in the above proof, fix either a plus or minus
sign. Then \be\label{www}e^{12}\pm e^{34},\q e^{13}\pm e^{42},\q e^{14}\pm
e^{23}\ee constitutes a basis of the 3-dimensional subspace $\E^2_\pm\bD$
giving rise to the celebrated decomposition
\be\label{SD}\E^2\bD=\E^2_+\bD\op\E^2_-\bD.\ee This determines a double
covering $SO(4)\to SO(3)_+\times SO(3)_-$, and there exist corresponding
\textit{subgroups} $SU(2)_+$, $SU(2)_-$ of $SO(4)$ acting trivially on
$\E^2_-\bD$, $\E^2_+\bD$ respectively.

The 2-form (\ref{hw}) belongs to the disjoint union \be\label{SS}\cS_+\sqcup
\cS_-, \ee where $\cS_\pm$ is a 2-sphere in $\E^2_\pm\bD$. The choice of sign
depends on whether $\hJ$ is positively or negatively oriented and is duplicated
in (\ref{56}). For example $J_0$ has fundamental 2-form \[
\ga=e^{12}+e^{34}+e^{56},\] and $\hg=e^{12}+e^{34}\in\cS_+$. The product
$\cS_+\times\cS_-$ may be identified with the Grassmannian of oriented 2-planes
in $\R^4$ and this was the origin of the concept of self-duality
\cite{AHS,ST}. Notice that (\ref{gs}) implies that $\hbox{Im}\,d$ lies in the
subspace $\E^2_+\bD$ of self-dual 2-forms; from this point of view $\bM$ is an
`instanton' over the torus $T^4$.\sb

The main result of \cite{AGS} may now be summarized by

\begin{theorem}\label{ags} The space $\cCgg$ is the disjoint union of
$\{J_0\}$ and the 2-sphere of all $g$-orthogonal almost complex structures $J$
on $\g$ for which $\hJ\in\cS_-$.\end{theorem}

Let \[\pZ=\{J\in\cCgg:\hJ\in\cS_-\}\] denote the 2-sphere featuring in this
theorem; we use the notation of \cite{AGS}. Consider $SO(4)$ as a subgroup of
$GL(6,\R)$ by letting it act trivially on $e^5,e^6$. Since $d(\g^*)$ is spanned
by 2-forms in $\E^2_+\bD$, the subgroup $SU(2)_-$ is a group of Lie algebra
automorphisms of $\g$, acting transitively on $\pZ$. This observation will be
important in \S4.

\setcounter{theorem}0
\subsection{Deformation of $J_0$}

The main purpose of what follows is to generalize Theorem~\ref{ags} by removing
the orthogonality constraint. We begin by decomposing the space of \textit{all}
almost complex structures on $\bD$ as \[\AC_+\sqcup\AC_-,\] where $\AC_\pm$
consists of those structures inducing a $\pm$ orientation on $\bD$. This is the
extension of (\ref{SS}) in the non-metric situation, and \[\AC_\pm\cong
\frac{GL^+(4,\R)}{GL(2,\C)}\supset \frac{SO(4)}{U(2)}\cong\cS_\pm.\] We then
set

\begin{definition} Let $\cC_\pm=\{J\in\cCpg:\hJ\in\AC_\pm\}$.
\end{definition}

\n In contrast to $\AC_\pm$, the definition of $\cC_\pm$ incorporates the
requirement of integrability. If the overall orientation of $\g$ is not fixed,
we obtain \[\cCg=\cC_+\sqcup\cC_-\sqcup(-\cC_+)\sqcup(-\cC_-),\] where
$-\cC_\pm=\{-J:J\in\cC_\pm\}$. Signs that appear as \textit{sub}scripts refer
exclusively to the orientation on $\bD$.

In order to gain a greater understanding of the subsets $\cC_+,\cC_-$, we now
describe a completely different set-theoretic partition of $\cCpg$, in which
$J_0$ plays the role of an origin. We use the notation (\ref{zz}), with
$\ow^{123}=\ow^1\we\ow^2\we\ow^3$ etc.

\begin{definition}\label{fin} Let\\[2pt] $\cCf_0$ be the open subset of 
$\cCpg$ consisting of complex structures admitting a basis $\{\a^1,\a^2,\a^3\}$
of $(1,0)$-forms for which $\a^{123}\we \ow^{123}\ne0$.\\[4pt] \n $\cCi_0$ be
the complement $\cCpg\sm\cCf_0$.\end{definition}

\n The zero subscript emphasizes that comparisons are being made with reference
to $J_0$, and elements of $\cCi_0$ are `infinitely far' from $J_0$ in the sense
that the coefficients in (\ref{uu}) below become unbounded.

\begin{proposition}\label{u} If $J\in\cCf_0$ then there exists a basis 
$\{\a^i\}$ of $(1,0)$-forms and $a,b,c,d,x,y$ in $\C$ such that
\be\label{uu}\left\{\ba{rcl}\a^1 &=& \w^1+a\ow^1+b\ow^2,\\\a^2
&=&\w^2+c\ow^1+d\ow^2,\\\a^3&=&\w^3+x\ow^1+y\ow^2+u\ow^3, \ea\right.\ee where
$u=-ad+bc$.\end{proposition}

\pr Theorem~\ref{JD} implies that $\a^1,\a^2$ can be chosen so that their real
and imaginary components span $\bD$. The condition $\a^{123}\we\ow^{123}\ne0$
ensures that $\w^1,\w^2,\w^3$ appear with non-zero coefficients. We thus obtain
the description (\ref{uu}) for some $u\in\C$. The equation relating $a,b,c,d,u$
is a direct consequence of the integrability condition \[d\a^3\we\a^{123}=0\]
expressing the fact that $d\a^3$ has no $(0,2)$-component.\qed

For the remainder of this section, we focus on $\cCf_0$. In (\ref{uu}), $\hJ$
is the almost complex structure on $\bD$ with $(1,0)$-forms
\be\label{ad}\left\{\ba{rcl}\a^1&=&\w^1+a\ow^1+b\ow^2,\\\a^2
&=&\w^2+c\ow^1+d\ow^2,\ea\right.\ee and is conveniently represented by the
matrix \be\X=\left(\!\ba{cc}a&b\\c&d\ea\!\right).\label{X}\ee The
characteristic polynomial of $\XoX$ has the form \[c(x)=x^2-\ga x+\delta,\]
where \[\ba{l}\ga= \tr(\XoX)=|a|^2+|d|^2+2\Re(b\ol
c),\\[3pt]\delta=\det(\XoX)=|u|^2.\ea\] Let $\la,\mu$ denote the roots of
$c(x)$. An inspection of the coefficients $\ga, \delta$ shows that $\la,\mu$
are either both real of the same sign or complex conjugates. The following
result is less obvious.

\begin{lemma}\label{coin} If $\la,\mu$ are real and non-positive, then
$\la=\mu\le0$.\end{lemma}

\pr Let $z=b\ol c\in\C$. The inequality $0\le|z|+\Re(z)$ implies that
\[ 0\le|bc|+\Re(b\ol c)\le|bc-ad|+|ad|+\Re(b\ol c).\]
Thus, \[ 0\le2|ad-bc|+|a|^2+|d|^2+2\Re(b\ol c)=2\sqrt\delta+\ga.\] If $\ga<0$
then the discriminant $\ga^2-4\delta$ is negative and the roots of $c(x)$ are
not real.\qed

A direct calculation reveals that \[\a^{12}\we\oa^{12}=(1-\ga+\delta)
\w^{12}\we\ow^{12}=c(1)\w^{12}\we\ow^{12}.\] Whence\smallbreak

\begin{proposition}\label{eig} 
\[\ba{rcl}\a^{12}\we\oa^{12} &=& 4(1-\la)(1-\mu)e^{1234},\\[3pt]
\a^{123}\we\oa^{123} &=& \!\!-8i(1-\la)(1-\mu)(1-\la\mu)e^{12\cdots6}.\ea\]
\end{proposition}

This proposition implies that the sign of $(1-\la)(1-\mu)$
corresponds to that of $\cC_\pm$. Moreover, $\la\mu\ge0$, so that\sb

(i) $J\in\cC_+\cap\cCf_0\ \Rightarrow\ 0\le \la\mu<1$,\sb

(ii) $J\in\cC_-\cap\cCf_0\ \Rightarrow\ \la\mu>1$.\sb

\n Note that $\mu=\ol\la$ implies that $c(1)=|1-\la|^2>0$, and is only
admissible for $J\in\cC_+$.

\vspace{-70pt}
\hspace{30pt}\scalebox{0.7}{\includegraphics{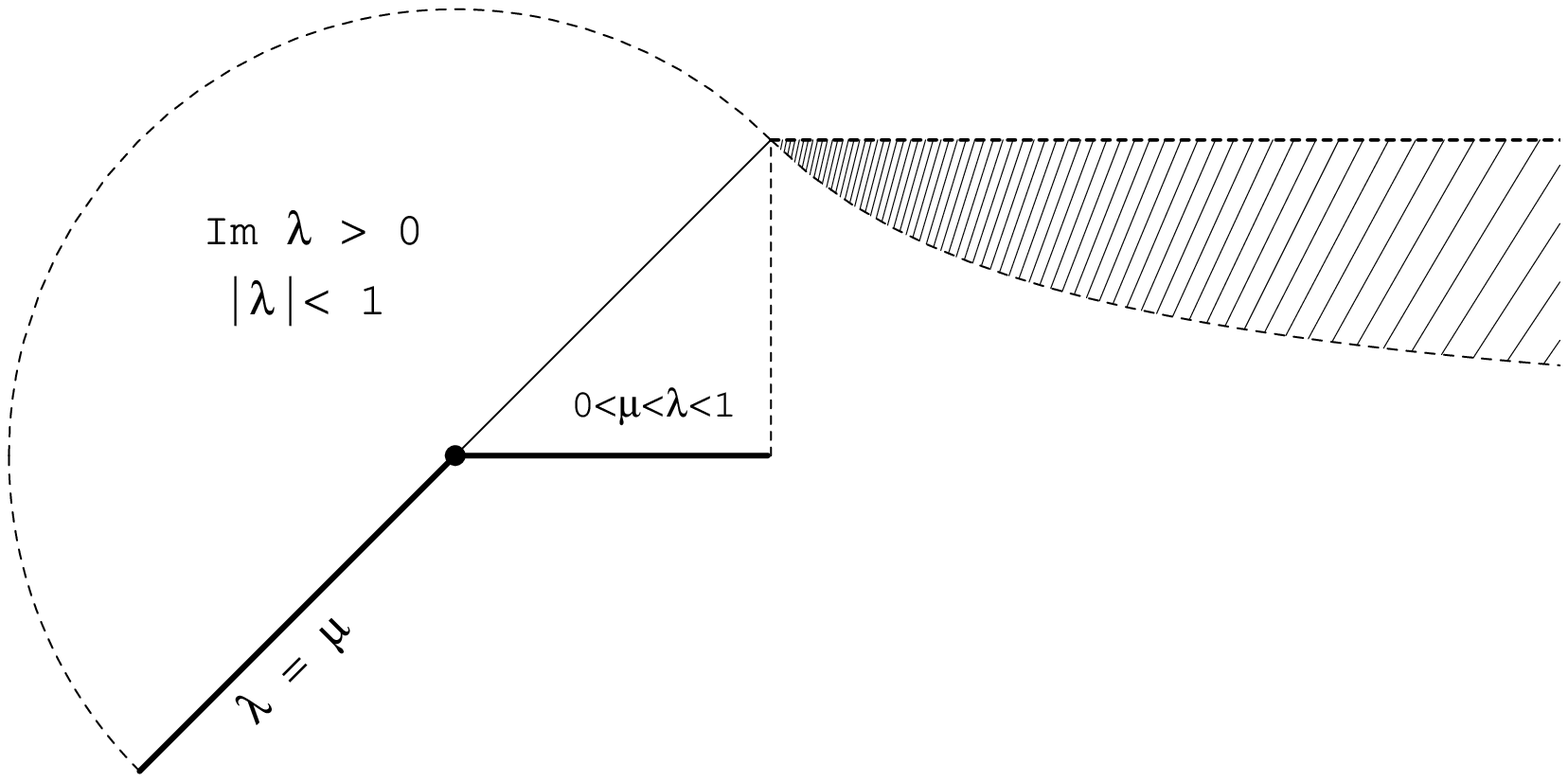}}
\vspace{-60pt}

\centerline{Figure~1}\vspace{20pt}

The possibilities for the unordered pair $\{\la,\mu\}$ are illustrated
schematically in Figure~1. The two labelled regions correspond to condition
(i), with the origin a common point of intersection. The semi-circular region
corresponds to $\hbox{Im}\,\la>0$ and $|\la|<1$. By contrast, points
$(\la,\mu)$ south-east of the diagonal line $\la=\mu$ represent those of the
real plane in the usual way: the triangular region bounds points arising from
$\cC_+$ with $\la>\mu\ge0$ and the shaded region represents points satisfying
(ii).\smallbreak

\re The similarity class of $\XoX$ is invariant by the action \[X\mapsto
g^{-1}X\ol g,\q g\in GL(2,\C),\] that characterizes the relation of
`consimilarity' \cite{HH}. This action is known to be transitive on the set of
$X$ corresponding to a fixed similarity class of $\XoX$, provided $X\ol X\ne0$.
In particular, if $\XoX$ is diagonalizable with $\la,\mu$ positive then there
exists $g\in GL(2,\C)$ such that \[X=g^{-1}\left(\!\ba{cc} \sqrt\la&0\\
0&\sqrt\mu \ea\!\right)\ol g\] (the sign of the square roots can be changed by
modifying $X$) \cite[4.6.11]{HH}. Points in the open triangular region
therefore represent $GL(2,\C)$ orbits of $\AC_+$ that consist of projections
(via Theorem~\ref{JD}) of complex structures in $\cC_+$. On the other hand,
$\cC_-$ contains elements for which $\la$ is infinite and the corresponding
eigenvector of $\XoX$ determines a point of $\mathbb{CP}^1$ that re-appears as
a 2-sphere in Theorem~\ref{2sph} below.\smallbreak

We are focussing attention on the set $\cCg$ of left-invariant complex
structures on $G$. The right action of $G$ induces a transitive action on $\bM$
and an induced action on $\cCpg$. Given an element $J$ of $\cCpg$, let
$R_{\G}(J)$ denote the orbit of $J$ induced by this action.

\begin{proposition} If $J$ is given by (\ref{uu}) and (\ref{X}), then
\[\dim_\C R_{\G}(J)=\left\{\ba{ll}0 & \hbox{if }X=0,\\
1&\hbox{if }u=0\hbox{ and }X\ne0,\\2&\hbox{otherwise}.\ea\right.\]
\end{proposition}

\pr Right translation leaves invariant the 1-forms $\w^1,\w^2$ in (\ref{zz}),
but maps $\w^3$ to $\w^3+p\w^1+q\w^2$ for arbitrary $p,q\in\C$. Thus
\[\ba{rcl}\a^3&\mapsto&\a^3+p\w^1+q\w^2+u(\ow^3+\ol p\ow^1+\ol q\ow^2)\\ &=&
\a^3+p(\a^1-a\ow^1-b\ow^2)+q(\a^2-c\ow^1-d\ow^2) +u(\ow^3+\ol p\ow^1+\ol
q\ow^2)\\ &=& \a^3+p\a^1+q\a^2+u\ow^3+(u\ol p-ap-qc)\ow^1+(u\ol
q-bp-dq)\ow^2.\ea\] This has the effect of replacing $(x,y)$ by $(x+u\ol
p-ap-cq,y+u\ol q-bp-dq)$ in (\ref{uu}). If $X=0$ then $J$ is unchanged, and
$R_{\G}(J)=\{J\}$. The remaining cases follow from the fact that $u=0$ if and
only if $ap+cq$ is proportional to $bp+dq$.\qed

A point of the moduli space of complex structures on $\bM$ consists of an
equivalence class of a complex structure (invariant or not) under the action of
the diffeomorphism group. A neigbourhood of it at a smooth point $J$ can be
identified with a subset of $H^1(\bM,\cO(T_J))$, where $T_J$ denotes the
holomorphic tangent bundle of $J$. This vector space is isomorphic to the
corresponding cohomology group of the Dolbeault complex \be\label{Dol}
0\to\W^{0,0}(T)\to\W^{0,1}(T)\to\W^{0,2}(T)\to\W^{0,3}(T)\to0.\ee Now, at
least if $J$ has rational coefficients relative to the basis $\{e^i\}$, it is
known that the cohomology of (\ref{Dol}) coincides with that of the
finite-dimensional subcomplex formed by restricting to left-invariant forms of
type $(p,q)$ \cite{CF,Cord}.

The cohomology of the invariant subcomplex is easily computed in the case of
$\bM$, using the techniques of \cite{S6}. In all cases, $\ker\ol\partial:
\W^{0,1}(T)\to\W^{0,2}(T)$ has dimension 6, whereas $\ol\partial(\W^{0,0}(T))$
has dimension 0,1,2, consistent with the above proposition. This phenomenon
leads to the jumping of Hodge numbers at $J_0$ described in \cite{Nak}. In any
case, it implies that the true moduli space of complex structures on $\bM$ has
dimension $4$ at generic points. It is also suggests that every point is
represented by an invariant complex structure, though the moduli space is
singular at $J_0$ and other boundary points in Figure 1.

\subsection{Study of $\cC_+$}

In recovering $J$ from $\hJ$, we need only worry about the coefficients of
$\a^3$ in (\ref{uu}), for which $u$ is determined by $a,b,c,d$ and $x,y$ are
arbitrary complex numbers. Connectivity properties of $\cCpg$ are determined by
those of its dense subset $\cCf_0$, and one expects the topology to be captured
by that of the domains (i),(ii) characterizing the choice of $\{\la,\mu\}$.
Results in this and the next sections will confirm that $\cC_+$ and $\cC_-$ are
the connected components of $\cCpg$.

\begin{proposition}\label{cap} $\cC_+\cap\cCi_0=\es$.\end{proposition}

\pr Let $J\in\cCi_0$. Suppose that $\{\a^1,\a^2,\a^3\}$ is a basis of
$(1,0)$-forms of $J$ with the real and imaginary components of $\a^1,\a^2$
spanning $\bD$. Consider the two cases:\sb

\n(i) $\a^{12}\we\ow^{12}\ne0$. This implies that $\a^3\in\al\ow^1,\ow^2,
\ow^3\ar$. The positive overall orientation of $J$ then forces $\hJ\in\AC_-$,
and $J$ cannot be in the same connected component as $J_0$.\sb

\n(ii) $\a^{12}\we\ow^{12}=0$ and $\a^3\not\in\al\ow^1,\ow^2,\ow^3\ar$. Then
$\al\a^1,\a^2\ar\cap\al\ow^1,\ow^2\ar\ne\{0\}$, and $\hJ$ has a non-zero
$(1,0)$-form $A\ow^1+B\ow^2$, which (without losing generality) we may take to
equal $\a^1$. If \[\a^2=P\w^1+Q\w^2+C\ow^1+D\ow^2,\] then the integrability of
$J$ forces $AD-BC=0$, and (subtracting a multiple of $\a^1$) we may suppose
that $C=D=0$. But then $\hJ\in\AC_-$, and again $J\not\in\cC_+$.\qed

\begin{theorem} $\cC_+$ is isomorphic to $\cU\times\C^2$ where 
$\cU$ is a star-shaped subset of $\C^4$.\end{theorem}

\pr Proposition~\ref{cap} implies that $\{\hJ:J\in\cC_+\}$ is a subset of
\[\cU=\{X\in\C^4:(1-\la)(1-\mu)>0,\ 0\le \la\mu<1\},\] using the notation 
(\ref{X}).  We shall show that if $X\in\cU$ then $tX\in\cU$ for any
$t\in[0,1]$, a fact that is illustrated by Figure~1. Indeed, if the eigenvalues
$\la,\mu$ of $\X\ol\X$ are complex conjugates, then the defining condition for
$\cU$ is $|\la|<1$, which becomes $t^{2}|\la|<1$ and remains valid. Suppose now
that $\la,\mu\in\R$. Then $(1-\la)(1-\mu)$ becomes
\[f=(1-t^{2}\la)(1-t^{2}\mu),\] an expression with roots $t_1^2=1/\la$ and
$t_2^2=1/\mu$. If $\la,\mu$ are both negative then $f$ has no real roots and is
always strictly positive. If $\la,\mu$ are both positive then at least one of
$1/\la,1/\mu$ is greater than 1, and $(1-\la)(1-\mu)>0$ implies that both are
greater than 1. It follows that $f>0$ for all $t\in[0,1]$, as required.

The restriction of $\p$ to $\cC_+$ is a trivial bundle, whose fibre is obtained
by varying only $x,y$, and $\cC_+$ can be identified with $\cU\times\C^2$.\qed

The complex structure induced on $\cCf_0$ and $\cU\times\C^2$ by the
coefficients in (\ref{uu}) obviously coincides with that induced by the natural
inclusion \[\cC(\g)\to\mathbb{G}\mathrm{r}_3(\C^6)\] obtained by mapping an
invariant complex structure $J$ to the span of a $(3,0)$-form $\a^{123}$. This
is also the natural complex structure induced from that of the potential
tangent space $H^1(\bM,\cO(T_J))$ to the moduli space \cite{S6}. From this
point of view, as a complex manifold, $\cC_+$ can be identified with an open
set of the quadric in $\C^7$ defined by the equation $u=-ad+bc$.

\re A completely different approach to describing complex structures on a
6-dimensional nilmanifold is based on properties of a $(3,0)$-form
$\a^{123}=\phi+i\psi$. The real component $\phi$ is a closed 3-form belonging
to the open orbit $\cO$ of elements of $\E^3\g^*\cong\R^{20}$ with stabilizer
isomorphic to $SL(3,\C)$. As a consequence, any element $\phi$ of $\cO$
determines a corresponding almost complex structure $J_\phi$ and
$\psi=J_\phi\phi$ \cite{Hit}. The kernel of $d:\E^3\g^*\to\E^4\g^*$ has
dimension 15, and $d(J_\phi\phi)=0$ turns out to be a single cubic equation in
the coefficients of $\phi$. This provides a description \[\cCg\cong\{\phi\in
\ker d\cap\cO: d(J_\phi\phi)=0\}/\C^*.\] More details will appear elsewhere.\mb

Consider an element $J\in\cC_+$ whose restriction to $\bD$ is $g$-orthogonal
and therefore an element of $\cS_+$. A point of $\cS_+$ at `finite' distance
from $\widehat{J_0}$ is given by (\ref{ad}) with $g(\a^i,\a^i)=0$ for $i=1,2$
and $g(\a^1,\a^2)=0$. This implies that $a=d=0$ and $b=-c$. It follows that the
space of $(1,0)$-forms of $J$ has a basis \be\label{PQ}\left\{\ba{rcl}
\b^1&=&\w^1+b\ow^2,\\\b^2 &=& -b\ow^1+\w^2,\\\b^3
&=&\w^3+x\ow^1+y\ow^2-b^2\ow^{3}.\ea\right.\ee 
Thus $\{\la,\mu\}=\{-|b|^2\}$, and $|b|<1$.\sb

\begin{corollary}\label{hemi} $\{\hJ:J\in\cC_+\}\cap\cS_+$ is an open 
hemisphere.\end{corollary}
\newcommand{\fq}[1]{\frac{#1}{1+|b|^2}} 

\pr From (\ref{PQ}), an element $\hJ\in\cS_+$ has $(1,0)$-forms
\[\left\{\ba{rcl}\a^1&=&e^1+ie^2+be^3-ibe^4,\\\a^2&=&\!\!-be^1+ibe^2+e^3+ie^4.
\ea\right.\] Setting \[A=\fq{1-|b|^2},\q B=i\fq{\ol b-b},\q C=-\fq{b+\ol b}\]
gives $A^2+B^2+C^2=1$ and \[\fq{\a^1-\ol b\a^2}=e^1+i(Ae^2+Be^3+Ce^4).\] In the
notation (\ref{www}) with plus signs, the fundamental 2-form of $\hJ$ equals
\[e^1\we(Ae^2+Be^3+Ce^4)+\cdots\ =A\w_1+B\w_2+C\w_3.\] The condition $|b|<1$
translates into $A>0$, that describes a hemisphere in $\cS_+$.\qed

\ex The almost complex structure $I$ on $\bD$ with space of $(1,0)$-forms \[\al
e^1+ie^3,\,e^4+ie^2\ar=\al \w^1+i\ow^2,\,\w^2-i\ow^1 \ar\] has $b=i$ in
(\ref{PQ}) and is a point on the equator $A=0$ of $\cS_+$. If $I$ were to equal
$\hJ$ with $J\in \cC_+$, then $J$ has a $(1,0)$-form of type \[\a^3=\w^3+x\w^1
+y\w^2+\ow^3,\] with the final coefficient $+1$ necessary to satisfy the
integrability condition. But then $\a^{12}\we\oa^{12}\we\a^3\we\oa^3=0$, which
is impossible.\mb

\setcounter{theorem}0
\subsection{Study of $\cC_-$}

We have remarked (Theorem~\ref{ags}) that the imposition of the standard metric
implies that $J_0$ is the only orthogonal structure in its component
$\cC_+$. Whilst $J_0$ is convenient for the study of $\cC_+$, it is less so for
$\cC_-$. For example, all the points of $\pZ$ belong to $\cCi_0$, making
calculations difficult in the coordinates of (\ref{uu}). We shall therefore
re-formulate Definition~\ref{fin} with respect to one particular element in
$\pZ$.

\begin{definition} Let\\[2pt] $J_1\in \cC_-$ denote the complex structure for
which $\eta=\w^1\we\ow^2\we\ow^3$ is a $(3,0)$-form,\\[4pt] $\cCf_1$ be the
open subset of $\cCpg$ consisting of complex structures admitting a basis
$\{\b^1,\b^2,\b^3\}$ of $(1,0)$-forms for which $\b^{123}\we\ol\eta \ne0$, and
\\[4pt]\n $\cCi_1$ be the complement $\cCpg\sm\cCf_1$.  \end{definition}

The analogue of Proposition~\ref{u} is

\begin{proposition}\label{bbb} If $J\in\cCf_1$ then $\b^i$ may be chosen so 
that \be\label{av}\left\{\ba{rcl}\b^1 &=& \w^1+a\ow^1+b\w^2,\\\b^2 &=&
\ow^2+c\ow^1+d\w^2,\\\b^3 &=& \ow^3+x\ow^1+y\w^2+v\w^3,\ea\right.\ee where
$a,b,c,d,x,y,v\in\C$ and $d=-av$.\end{proposition}

\pr This follows from Theorem~\ref{JD} and the equation $d\b^3\we\b^{123}=0$.
\qed

Because of the equation $d=-av$, $v$ is unconstrained if $a$ happens to vanish,
and this contrasts with the situation in the previous section. Let $\ACf_1$ be
the set of almost complex structures on $\bD$ with a basis of $(1,0)$-forms
consisting of \be\label{abcd}\left\{\ba{rcl}\b^1 &=&\w^1+a\ow^1+b\w^2,\\[4pt]
\b^2 &=&\ow^2+c\ow^1+d\w^2,\ea\right.\ee for some $a,b,c,d\in\C$. The almost
complex structure $\widehat{J_1}$ corresponds to $a=b=c=d=0$.

\ex Recall that the projection $J\mapsto\hJ$ maps $\pZ$ onto $\cS_-$. Elements
of $\cS_-$ have the form (\ref{abcd}) with $a=d=0$ and $b=-c$, except that
$-\widehat{J_1}$ corresponds to $b$ and $c$ infinite. Thus, any element of
$\cS_-\sm\{-\widehat{J_1}\}$ equals $\hJ$ for some $J\in\cCf_1$ (compare 
Corollary~\ref{hemi}).\mb

With respect to the basis $\{e^1,e^2,e^3,e^4\}$ of $\bD$, the element
$\widehat{J_1}$ is represented by the matrix \[Q_1=
\left(\ba{cccc}0&-1&0&0\\1&0&0&0\\0&0&0&1\\0&0&-1&0\ea\right)\in SO(4).\] We
can then identify $\AC_-$ with the orbit $\{X^{-1}Q_1 X:X\in GL^+(4,\R)\}$, any
element of which admits a polar decomposition \be\label{pd}X^{-1}Q_1X=SP,\ee 
where $S$ is symmetric positive-definite and $P\in SO(4)$.

\begin{lemma}\label{retr} With the above notation, $P^2=-1$, and 
the resulting mapping $\ret\colon\AC_-\to\cS_-$ defined by $SP\mapsto P$ is a
retraction.\end{lemma}

\pr By first diagonalizing $S$, we may find a symmetric matrix $\si$ for which
\[S=e^\si=\sum_{k=0}^\infty\frac1{k!}\si^k.\] We claim that $e^{t\si}P\in
\AC_-$ for all $t\in[0,1]$. Since $(SP)^2=-1$, \[SP=-P^{-1}S^{-1}=(P^{-1}
S^{-1}P) (-P^{-1}),\] in which $P^{-1}S^{-1}P=P^{T} S^{-1}P$ is
positive-definite symmetric. Uniqueness of the polar decomposition implies that
$S=P^{-1}S^{-1}P$ and $P=-P^{-1}$, so $P^2=-1$. It follows also that
$\si=-P^{-1}\si P$ and \[e^{t\si}P=(P^{-1}e^{-t\si}P)P= -P^{-1}(e^{t\si})^{-1}
=-(e^{t\si}P)^{-1},\] as required.\qed

\begin{proposition} $\ret^{-1}(Q_1)\cap\{\hJ:J\in\cCi_1\}=\es$.
\end{proposition}

\pr Given an element $SQ_1$ of $\ret^{-1}(Q_1)\cap\ACf_1$ with $(1,0)$-forms as
in (\ref{abcd}), we claim that $b=c$. Identifying almost complex structures
with $4\times4$ matrices, $1+iSQ_1$ annihilates the $(1,0)$-forms $\b^1,\b^2$
of (\ref{abcd}). Extending the standard metric $g$ on $\bD$ to a complex
bilinear form, \[\ba{rcl}0 &=& g((1+iSQ_1)\b^1,Q_1\b^2)
-g(Q_1\b^1,(1+iSQ_1)\b^2)\\[2pt] &=& 2g(\b^1,Q_1\b^2)\\[2pt]&=& 2i(b-c),\ea\]
as stated.

In analogy to Proposition~\ref{eig}, we have \be\b^{12}\we\ol\b^{12}=-(1-\la)
(1-\mu)e^{1234},\label{bb}\ee where $\la,\mu$ are the eigenvalues of $Y\ol Y$
with $Y=\hbox{\small$\left(\!\ba{cc}a&b\\b&d\ea\!\right)$}$. Since $Y\ol Y$ is
Hermitian, $\la,\mu$ are non-negative. Equation (\ref{bb}) forces $\la,\mu$ to
lie in the interval $[0,1)$, and $Y$ is bounded. Thus, $\ret^{-1}(Q_1)\subset
\ACf_1$. Finally suppose that $SQ_1=\hJ$ with $J\in\cCi_1$. From (\ref{av}),
it must be the case that $J$ has a $(1,0)$-form $\b^3$ belonging to the span of
$\ow^1,\w^2,\w^3$. But this is impossible, given that $J\in\cCpg$.\qed

\begin{theorem}\label{2sph} $\cC_-$ has the homotopy type of a
2-sphere. \end{theorem}

\pr Let $\cV=\{J\in\cC_-:\hJ\in\ret^{-1}(Q_1)\}$. We first show that this space
is contractible. As a consequence of the previous proposition, a complex
structure $J$ in $\cV$ has a basis of (1,0)-forms
\[\label{fibres}\left\{\ba{rcl} \b^1 &=& \w^{1}+ta\ow^{1}+tb\w^{2}\\
\b^2 &=& \ow^{2}+tb\ow^{1}-t^2av\w^{2}\\\b^3 &=&
\ow^{3}+tx\ow^{1}+ty\w^{2}+tv\w^{3},\ea\right.\] for some $a,b,x,y,v\in\C$ and
(for the moment) $t=1$. But if we now allow $t$ to vary in the interval
$[0,1]$, these forms define a complex structure in $\cV$. This process defines
a homotopy \[\cV\times[0,1]\to\cV,\] with the property that $(J,0)$ maps to
$J_1$ for all $J$ in $\cV$.

The fact that all elements of $\pZ$ are equivalent under a $SU(2)_-$ action
(see the remarks at the end of \S1) allows us to extend the above to all the
fibres of $\ret$ over all points of the 2-sphere $\pZ$.\qed

Combined with the theorem in \S3, we can conclude

\begin{theorem} $\cCpg$ has the same homotopy type as $\cCgg$, where $g$ is 
the inner product (\ref{g}).\end{theorem}

\bigbreak\n\textit{Acknowledgments.} Thanks are due to John Rawnsley for
indicating the proof of Lemma~\ref{retr}, and for help in connection with the
remark in \S2.

\vskip10pt\small

\n Mathematical Institute, University of Oxford, 24--29 St Giles', Oxford, 
OX1~3LB, UK\vskip8pt

\n Dipartimento di Matematica, Politecnico di Torino, Corso Duca degli 
Abruzzi 24, 10129 Torino, Italy

\enddocument